\documentclass[12pt,leqno]{amsart}

\DeclareMathOperator{\cd}{cd}%
\DeclareMathOperator{\cor}{cor}%
\DeclareMathOperator{\rank}{rank}%
\DeclareMathOperator{\res}{res}%

\usepackage{amssymb}
\usepackage{hyperref}
\usepackage[arrow]{xypic}

\newcommand{\ann}{\text{\rm ann}}
\newcommand{\F}{\mathbb{F}}
\newcommand{\Fp}{\F_p}
\newcommand{\Gal}{\text{\rm Gal}}
\newcommand{\Ic}{\mathcal{I}}
\newcommand{\Q}{\mathbb{Q}}
\newcommand{\Sc}{\mathcal{S}}
\newcommand{\Tc}{\mathcal{T}}
\newcommand{\Z}{\mathbb{Z}}

\begin{document}

\title[Demu\v{s}kin Groups]
{Demu\v{s}kin groups, Galois modules, and the elementary type
conjecture}

\author[Labute]{John Labute}
\address{Department of Mathematics and Statistics, \
McGill University, \linebreak Burnside Hall, 805 Sherbrooke Street
West, \ Montreal, Quebec \linebreak H3A 2K6 \ CANADA}
\email{labute@math.mcgill.ca}

\author[Lemire]{Nicole Lemire$^{\dag}$}
\address{Department of Mathematics, Middlesex College, \
University of Western Ontario, London, Ontario \ N6A 5B7 \ CANADA}
\thanks{$^\dag$Research supported in part by NSERC grant R3276A01.}
\email{nlemire@uwo.ca}

\author[Min\'{a}\v{c}]{J\'an Min\'a\v{c}$^{\ddag}$}
\thanks{$^\ddag$Research supported in part by NSERC grant R0370A01,
by the Mathematical Sciences Research Institute, Berkeley, by the
Institute for Advanced Study, Princeton, and by a 2004/2005
Distinguished Research Professorship at the University of Western
Ontario.} \email{minac@uwo.ca}

\author[Swallow]{John Swallow}
\address{Department of Mathematics, Davidson College, Box 7046,
Davidson, North Carolina \ 28035-7046 \ USA}
\email{joswallow@davidson.edu}


\begin{abstract}
    Let $p$ be a prime and $F(p)$ the maximal $p$-extension of a
    field $F$ containing a primitive $p$-th root of unity.  We
    give a new characterization of Demu\v{s}kin groups among
    Galois groups $\Gal(F(p)/F)$ when $p=2$, and, assuming the
    Elementary Type Conjecture, when $p>2$ as well.  This
    characterization is in terms of the structure, as Galois
    modules, of the Galois cohomology of index $p$ subgroups of
    $\Gal(F(p)/F)$.
\end{abstract}

\date{December 12, 2005}

\maketitle

\newtheorem{theorem}{Theorem}
\newtheorem{proposition}{Proposition}
\newtheorem{lemma}{Lemma}
\newtheorem*{corollary*}{Corollary}
\newtheorem*{ETC*}{Elementary Type Conjecture for Odd $p$}

\theoremstyle{definition}
\newtheorem{remark}{Remark}

\parskip=10pt plus 2pt minus 2pt

Let $p$ be a prime and $F$ a field containing a primitive $p$-th
root of unity $\xi_p$.  The union $F(p)$ of all finite Galois
extensions $L/F$ in a fixed algebraic closure of $F$ with $[L:F]$ a
power of $p$ is called the maximal $p$-extension of $F$.  Consider
$G=\Gal(F(p)/F)$.  Observe that while every profinite group is a
Galois group of some Galois extension \cite{W1}, the condition that
$G=\Gal(F(p)/F)$ is substantially more restrictive.

We ask when $G=\Gal(F(p)/F)$ is a Demu\v{s}kin group, that is, a
finitely generated pro-$p$-group satisfying $\dim_{\Fp} H^2(G,\Fp) =
1$ such that the cup product
\begin{equation*}
    \gamma_F : H^1(G,\Fp) \times H^1(G,\Fp) \to H^2(G,\Fp)
\end{equation*}
is a nondegenerate bilinear form.  (See \cite[\S III.9]{NSW} and
\cite[\S I.4.5 and \S II.5.6]{S2}.  Others relax the requirement on
finite generation, but in this article we consider only the finitely
generated case.)  Demu\v{s}kin groups arise, for instance, as
pro-$p$-completions of fundamental groups of compact surfaces $T$ of
genus $g \ge 1$ when $T$ is orientable and $g \ge 2$ when $T$ is not
orientable.  If $F$ is a finite extension of the field of $p$-adic
numbers $\Q_p$ and contains $\xi_p$, then $G$ is Demu\v{s}kin.

The study of Demu\v{s}kin groups among Galois groups $\Gal(F(p)/F)$ is
an important part of the program of classification of possible Galois
groups of maximal $p$-extensions of fields, as these groups form an
essential part of the local theory of this project. In turn, the
classification of possible $\Gal(F(p)/F)$ is one of the key problems
in current Galois theory. This study is also crucial for the
development of anabelian algebraic geometry over fields. (See
\cite{Ef3}, \cite{Ko}, and further references in these papers.)

In this paper we detect whether $G=\Gal(F(p)/F)$ is a Demu\v{s}kin
group in terms of the Galois module structure of the Galois
cohomology of index $p$ subgroups of $G$. For such $G$, we establish
a new characterization of Demu\v{s}kin groups when $p=2$. When $p>2$
our characterization depends upon the Elementary Type Conjecture in
the theory of Galois pro-$p$-groups.  The close relationship of this
characterization with the Elementary Type Conjecture offers a new
approach to the Conjecture. (See Remark~\ref{re:etc} in
section~\ref{se:proofs}.)

The surprising new insight contained in Theorem~\ref{th:new} below
is that $G$ is Demu\v{s}kin if $\cd (G) = 2$ and $H^2(N,\Fp)$,
with $N$ a subgroup of index $p$ of $G$, does not grow ``too fast."
In fact a relatively mild condition on the growth of $H^2(N,\Fp)$
guarantees that $\dim_{\Fp} H^2(N,\Fp) = 1$.

In considering Galois cohomology groups as Galois modules, we could
use the results of \cite{LMS1} and \cite{LMS2}.  These results,
however, depend upon recent, complex, partially published work of
Rost-Voevodsky on the Bloch-Kato Conjecture. For the proof of the
following theorem we use only the results in \cite{MeSu} concerning
the Bloch-Kato Conjecture in the case $n=2$.

Before formulating the main theorem we recall that if $G$ is a
finitely generated pro-$p$-group, any subgroup of index $p$ is closed
\cite[\S I.4.2, Ex. 6]{S2}, and that in a pro-$p$-group, any subgroup
of index $p$ is normal. Let $H$ be a group and $M$ be an
$\Fp[H]$-module. We say that $M$ is a trivial $\Fp[H]$-module if for
each $\tau \in H$ and $m \in M$ we have $\tau(m) = m$.
\begin{theorem}\label{th:new}
    Let $F$ be a field containing a primitive $p$-th root of unity,
    and suppose that $G=\Gal(F(p)/F)$ is a finitely generated
    pro-$p$-group of cohomological dimension $2$. Then for each
    subgroup $N$ of $G$ of index $p$, the following conditions on
    the $\Fp[G/N]$-module $H^2(N,\Fp)$ are equivalent:
    \begin{enumerate}
        \item $H^2(N,\Fp)$ has no nonzero free summand.
        \item $H^2(N,\Fp)$ is a trivial $\Fp[G/N]$-module.
    \end{enumerate}

    Now assume additionally that either $p=2$ or $p>2$ and the
    Elementary Type Conjecture holds. (See the end of
    section~\ref{se:pquatetc}.)

    Then $G$ is Demu\v{s}kin if and only if, for every subgroup $N$
    of $G$ of index $p$, $H^2(N,\Fp)$ has no nonzero free summand.
\end{theorem}

We observe that the Elementary Type Conjecture has been established
for some important classes of fields, including algebraic extensions
$F$ of $\Q$ with finitely generated $G=\Gal(F(p)/F)$. (See
\cite{Ef1} and \cite{Ef2}.) For such fields Theorem~\ref{th:new} is
a precise characterization. For additional information about the
Elementary Type Conjecture see \cite{Ma2}.

In \cite[Theorem~1]{DuLa} it was shown that a finitely generated
pro-$p$-group $G$ such that $\dim_{\Fp} H^1(G, \Fp) > 1$ and
$\dim_{\Fp} H^2(G,\Fp) = 1$ is Demu\v{s}kin if and only if
$H^2(N,\Fp) \cong \Fp$ for all subgroups $N$ of $G$ of index $p$.
Note that in Theorem~\ref{th:new}, under the assumption that
$\cd(G)=2$ and, if $p>2$, that the Elementary Type Conjecture holds,
we do not require that $H^2(N,\Fp)\cong \Fp$ but instead only that
$H^2(N,\Fp)$ contains no nonzero free summand. In fact, we prove
more than we claim in Theorem~\ref{th:new}. Namely, from the proof
of Theorem~\ref{th:new} it follows that we can replace the
hypothesis $\cd(G) = 2$ by two conditions which follow from it:
first, that the corestriction map from $H^2(N,\Fp)$ to $H^2(G,\Fp)$
is surjective for all subgroups $N$ of $G$ of index $p$, and,
second, that $H^2(G,\Fp)$ is not zero. We use deep results from
Galois cohomology in our proof, and it would be interesting to see
whether these characterizations of Demu\v{s}kin groups among groups
$\Gal(F(p)/F)$ also hold in the category of pro-$p$-groups.

The heart of our analysis is section~\ref{se:sc}, where we determine
the structure of the $\Fp[G/N]$-module $H^2(N,\Fp)$ when $N$ is a
subgroup of $G$ of index $p$ and the corestriction map $\cor:
H^2(N,\Fp)\to H^2(G,\Fp)$ is surjective.  In particular, we show
that $H^2(N,\Fp)\cong X\oplus Y$ where $X$ is trivial and $Y$ is a
free $\Fp[G/N]$-module.  Moreover, we characterize all such
decompositions of $H^2(N,\Fp)$.  We believe that these results are
of independent interest; for example, we obtain immediately from
this structure some information on the size of $H^2(N,\Fp)$.  (See
the Corollary to Theorem~\ref{th:h2}.)

Our approach uses $p$-quaternionic pairings, and we closely follow
\cite{Ku} in the first two sections. Some of the basic concepts
recalled here were introduced by Hwang and Jacob in \cite{HJ}.  In the
next sections we consider $H^2(G,\F_2)$ when cup products are strongly
regular, as well as the Galois module structure of $H^2(N,\Fp)$ for
index $p$ subgroups $N$ of $G$ when the corestriction is surjective.
Then we prove Theorem~\ref{th:new} and its Corollary.  Finally, we
close with a consideration of the $\Fp[G/N]$-module structure of the
1-cohomology groups $H^1(N,\Fp)$.

\section{$p$-Quaternionic Pairings}\label{se:pquat}

We seek to understand the condition on the cup product in the
definition of Demu\v{s}kin groups by considering such products in
the context of $p$-quaternionic pairings and bilinear forms in
general.

Let $H$ and $Q$ be elementary abelian $p$-groups written
multiplicatively and additively, respectively, and if $p=2$ choose
a distinguished element $-1\in H$, which may be trivial. Let
$\gamma: H\times H \to Q$ be a bilinear form. For a given element
$a \in H$, we define the group homomorphism
\begin{equation*}
    \gamma_a:H \to Q, \quad\quad \gamma_a(x) := \gamma(a,x), \ \ x
    \in H,
\end{equation*}
and we denote by $Q(a)$ the value group $\gamma_a(H)$ of
$\gamma_a$. We also define
\begin{equation*}
    N(a) = N_\gamma (a) = \ker \gamma_a = \left\{b \in H \ \ \vert\
    \ \gamma(a,b) = 0 \right\}.
\end{equation*}
We have $H/N(a) \cong Q(a)$.

The bilinear form $\gamma$ is called \emph{nondegenerate} if $Q(a)
\neq \{0\}$ for all $a\in H\setminus \{1\}$. If $Q(a) = Q$ for all
$a \in H \setminus \{1\}$, then the bilinear mapping $\gamma$ is
called \emph{strongly regular}.

Observe that in the following definition of $p$-quaternionic
pairing, the distinguished involution $-1 \in H$ is necessarily
$1$ if $p > 2$.

We say that $(H,Q,\gamma)$ is a $p$-\emph{quaternionic pairing}
if there exists a distinguished involution $-1 \in H$ such that:
\begin{enumerate}
    \item $Q$ is generated by the union $\bigcup_{a\in H}Q(a)$ of
    the value groups; \item $\gamma(a,a) = \gamma(a,-1)$ for all $a
    \in H$; \item if $p=2$ then $\gamma$ satisfies \emph{the linkage
    condition}: if for $a, b, c, d\in H$ we have that $\gamma(a,b) =
    \gamma(c,d)$, then there exists $e \in H$ such that
    \begin{equation*}
        \gamma(a,b) = \gamma(a,e) = \gamma(c,e) = \gamma(c,d);
    \end{equation*}
    and
    \item for every $n\ge 2$, the \emph{$M(n)$ condition} holds: if
    for elements $a_1, a_2, \dots, a_n$ and $b_1, b_2, \dots, b_n$
    of $H$ with $a_1,\dots,a_n$ linearly independent over $\Fp$, we
    have
    \begin{equation*}
        \sum_{i=1}^{n}\gamma(a_i,b_i) = 0,
    \end{equation*}
    then there exist elements $a_{n+1},\dots,a_k \in H$, with
    $a_1,a_2,\dots,a_k$ linearly independent over $\Fp$ and
    $x_{j_1,\dots,j_k} \in N_\gamma \left(a_1^{j_1} a_2^{j_2} \cdots
    a_k^{j_k}\right)$ such that
    \begin{equation*}
        b_i = \prod_{0 \le j_1, j_2, \dots, j_k \le
        p-1} \left( x_{j_1,\dots,j_k} \right)^{j_i},
    \quad i=1, \dots, k,
    \end{equation*}
    where $b_{n+1},\dots, b_k = 1$.
\end{enumerate}

(It is worth observing that from the bilinearity of $\gamma$
and condition (2), it follows that $\gamma$ is skew-symmetric
if $p>2$ and symmetric if $p=2$.)

A $p$-quaternionic pairing $(H,Q,\gamma)$ is said to be
\emph{strongly regular} if $\gamma$ is strongly regular. A
$p$-quaternionic pairing is said to be \emph{finite} if $H$ is
finite.

We consider several types of $p$-quaternionic pairings.

\emph{The cup product $\gamma_F$ of a field $F$}. Let $F$ denote a
field containing $\xi_p$ and let $G=\Gal(F(p)/F)$. The cup product
pairing
\begin{equation*}
    \gamma_F: H^1(G,\Fp) \times H^1(G,\Fp) \to H^2(G,\Fp)
\end{equation*}
satisfies the $M(n)$ conditions (see \cite[Proposition~4]{Me} and
\cite[(11.5) Theorem]{MeSu}). In fact, the $M(n)$ conditions are a
translation of the condition for the splitting of a sum of $n$
symbols in Milnor's $k_2 F=K_2 F/p K_2 F$ to the language of
$p$-quaternionic pairings.  Observe that in \cite{Ku} the set of
conditions $M(n), n \ge 2$, is a different condition than our
condition~(4), but the alteration of axiom (3) in \cite[p.~40]{Ku}
does not affect the results from \cite{Ku} that we use.
Condition~(1) also follows from \cite[(11.5)~Theorem]{MeSu}.
Condition~(2) is true when the distinguished involution $(-1) \in
H^1 (G,\Fp)$ corresponds to $[-1] \in F^\times/F^{\times p}$ via
Kummer theory. The linkage condition in (3) is well-known. (See
\cite[Chapter~3, Theorem~4.13]{L}.)

\emph{Pairings of $p$-local type}. A finite $p$-quaternionic pairing
$(H,Q,\gamma)$ with $Q$ a group of order $p$ and $\gamma$
nondegenerate is said to be \emph{of $p$-local type}.  If
$G=\Gal(F(p)/F)$ is a Demu\v{s}kin group, then it follows from the
definition that $(H^1(G,\Fp), H^2(G,\Fp), \gamma_F)$ is a
$p$-quaternionic pairing of $p$-local type.

\emph{Strongly regular pairings}. Suppose that $(H,\Fp,\gamma)$ is
a $p$-quaternionic non-degenerate pairing. Observe that for each
$a\in H\setminus \{1\}$, the subgroup $N(a)$ of $H$ is of index
$p$. Moreover, for $a,b \in H \setminus \{1\}$ we have
\begin{equation*}
    N(a) N(b) \neq H \Longleftrightarrow N(a) = N(b)
    \Longleftrightarrow \langle a \rangle = \langle b
    \rangle.
    \end{equation*}
Hence such $p$-quaternionic pairings are strongly regular. Now
suppose instead that $Q = \{0\}$. If we set $\gamma (a,b) = 0$ for
all $a,b \in H$, then $(H,Q,\gamma)$ is a $p$-quaternionic
pairing, called \emph{totally degenerate}. Totally degenerate
pairings are also strongly regular.

\emph{Pairings of weakly $p$-local type}. A $p$-quaternionic pairing
with $H = \{1\}$ is called \emph{trivial}. Each trivial pairing is
totally degenerate. We say that totally degenerate $p$-quaternionic
pairings, as well as pairings of $p$-local type, are pairings
\emph{of weakly $p$-local type}.

\section{$p$-Quaternionic Pairings and the Elementary Type
Conjecture}\label{se:pquatetc}

For $p>2$ we define the direct product and the group extension of
$p$-quaternionic pairings and consider the Elementary Type
Conjecture. Because we do not need it in Theorem~\ref{th:new}, we do
not consider the Elementary Type Conjecture when $p=2$. (See
\cite[Chap.~5]{Ma1} for the $p=2$ case in the context of abstract
Witt rings.)

(A) \emph{Direct product}.  Let $(H_i,Q_i,\gamma_i)$, $i=1, 2$, be
$p$-quaternionic pairings. Define $H=H_1 \times H_2$, $Q = Q_1
\times Q_2$, and
\begin{equation*}
    \gamma\big([a_1,a_2],[b_1,b_2]\big) = \big[\gamma_1(a_1,b_1),
    \gamma_2(a_2,b_2)\big], \quad a_i, b_i\in H_i.
\end{equation*}
Then $(H,Q,\gamma)$ is a $p$-quaternionic pairing called the
\emph{direct product}.

(B) \emph{Group extension}.  Suppose that $(H',Q',\gamma')$ is a
$p$-quaternionic pairing and let $T$ be a nontrivial finite
elementary abelian $p$-group. The \emph{group extension} of
$(H',Q',\gamma')$ by $T$ is the $p$-quaternionic pairing
$(H,Q,\gamma)$, where $H = H' \times T$, $Q = Q' \times (H'
\otimes T) \times (T \wedge T)$, and the pairing $\gamma:H
\times H \to Q$ is given by
\begin{equation*}
    \gamma \big([a_1,t_1],[a_2,t_2]\big) = \big[\gamma'(a_1,a_2),\ \
    a_1 \otimes t_2 - a_2 \otimes t_1,\ \ t_1 \wedge t_2\big].
\end{equation*}
Here $\otimes$ denotes the tensor product over $\Fp$ and $\wedge$
the exterior product.

For $p>2$, we say a that a finite $p$-quaternionic pairing is
\emph{of elementary type} if it may be constructed from
$p$-quaternionic pairings of weakly $p$-local type using the
operations of (a)~direct product and (b)~group extension by
nontrivial elementary abelian $p$-groups. The Elementary Type
Conjecture for $p>2$ is then as follows.  (We note that there are
several variants of the Elementary Type Conjecture which aim at the
classification of finitely generated $\Gal(F(p)/F)$, contained in
\cite{Ef1}, \cite{Ef2}, \cite{En}, \cite{JW}, and
\cite[p.~123]{Ma1}.)
\begin{ETC*}
    Let $p>2$ be a prime and $F$ a field containing a primitive
    $p$-th root of unity.  Suppose that $G=\Gal(F(p)/F)$ is a
    finitely generated pro-$p$-group. Then the cup product pairing
    $\gamma_F$ is of elementary type.
\end{ETC*}

\begin{theorem}[{\cite[Corollary~5]{Ku}}]\label{th:notstrreg}
    For $p>2$, a $p$-quaternionic pairing of elementary type is
    not strongly regular unless it is of weakly $p$-local type.
\end{theorem}

\section{Strongly regular cup products and $H^2(G,\F_2)$}
For the proof of the following proposition we originally used
streamlined arguments from \cite[pp.~42--43]{FY}.  Afterwards Kula
sent us a nice simplification of the proof, using ideas in
\cite[Proof of Proposition~2.16]{K1}.  We are grateful to him for
permitting us to adapt this simplification for use here.

\begin{proposition}\label{pr:p2}
    Let $F$ be a field of characteristic not $2$, and suppose that
    $G = \Gal(F(2)/F)\neq \{1\}$ is a finitely generated
    pro-$2$-group with $\gamma_F$ nondegenerate and strongly regular.
    Then $H^2(G,\F_2)\cong \F_2$.
\end{proposition}

\begin{proof}
    Assume that the hypotheses of our proposition are valid, and
    denote by $\vert A\vert$ the cardinality of a set $A$. Because
    the statement is trivial in the case $\vert H^1(G,\F_2)\vert
    =2$, we assume without loss of generality that $g := \vert
    H^1(G,\F_2)\vert >2$.  Denote $h:=\vert H^2(G,\F_2)\vert>1$, as
    $\gamma_F$ is nondegenerate.  Set $\ann (a)=\{(b)\in
    H^1(G,\F_2)\ \vert \ (a)\cdot (b)=0\}$. Since $(a)\cdot
    H^1(G,\F_2)\cong H^1(G,\F_2)/\ann (a)$, we see that
    \begin{equation*}
        \vert\ann(a)\vert = \frac{\vert H^1(G,\F_2)\vert}{\vert
        (a)\cdot H^1(G,\F_2)\vert} = \frac{\vert
        H^1(G,\F_2)\vert}{\vert H^2(G,\F_2)\vert} = \frac{g}{h}
    \end{equation*}
    for all nonzero $(a)\in H^1(G,\F_2)$.

    We show now that for arbitrary distinct, nonzero elements $(a),
    (b)\in H^1(G,\F_2)$, we have $\ann (a)+\ann (b)=H^1(G,\F_2)$.
    Let $(x)\in H^1(G,\F_2)$ be arbitrary.  If $q:=(x)\cdot (a)=0$,
    then $(x)\in \ann (a)$.  Assume therefore that $q\neq 0$.  Using
    the surjectivity of the map $((a)+(b)) \cdot -: H^1(G,\F_2)\to
    H^2(G,\F_2)$ and the linkage property, we see that there exists
    $(c)\in H^1(G,\F_2)$ such that
    \begin{equation*}
        q = (a)\cdot (x) = (a)\cdot (c) = ((a)+(b))\cdot (c).
    \end{equation*}
    Hence $((x)+(c))\cdot (a)=0=(b)\cdot(c)$ and therefore
    $(x)+(c)\in \ann (a)$ and $(c)\in \ann (b)$. Thus
    $(x)=((x)+(c))+(c)\in \ann(a) +\ann(b)$, as required.

    Let $D$ be the set of nonzero elements in the dual space of
    $H^1(G,\F_2)$.  Similarly, for each nonzero element $(a)\in
    H^1(G,\F_2)$, let $D_{(a)}$ be the set of all maps in $D$
    which are zero on $\ann (a)$.  Because $\ann (a) + \ann (b) =
    H^1(G,\F_2)$ for all pairs of distinct, nonzero elements $(a)$
    and $(b)$, $D$ contains the disjoint union of all $D_{(a)}$ in
    $D$. Since $\vert D\vert = g-1$ and $\vert D_{(a)}\vert=h-1$ for
    each nonzero $(a)$, we obtain $(g-1)(h-1)\le (g-1)$.  Therefore
    $h=2$.
\end{proof}

\section{Surjective Corestrictions and $H^2(N,\Fp)$}\label{se:sc}

In the following theorem we do not assume that $G$ is finitely
generated.

\begin{theorem}\label{th:h2}
    Let $F$ be a field containing a primitive $p$-th root of unity,
    and suppose that $G=\Gal(F(p)/F)$.  Let $N$ be a subgroup of $G$
    of index $p$, and suppose that the corestriction map
    $\cor:H^2(N,\Fp)\to H^2(G,\Fp)$ is surjective. Let $a\in
    F^\times$ be chosen so that the fixed field of $N$ is
    $F(\root{p}\of{a})$. Then the $\Fp[G/N]$-module $H^2(N,\Fp)$
    decomposes as
    \begin{equation*}
        H^2(N,\Fp) = X \oplus Y
    \end{equation*}
    where $X$ is a trivial $\Fp[G/N]$-module, $Y$ is a free
    $\Fp[G/N]$-module, and
    \begin{enumerate}
        \item $\dim_{\Fp} X = \dim_{\Fp} H^1(G,\Fp)/
        \ann(a)$
        \item $\rank_{\Fp[G/N]} Y = \dim_{\Fp} H^2(G,\Fp)/
        (a)\cdot H^1(G,\Fp)$.
    \end{enumerate}
\end{theorem}
After the proof, we characterize in Theorem~\ref{th:h2bis} all
decompositions of $H^2(N,\Fp)$ into direct sums of trivial
and free submodules.

Observe that we have a natural sequence
\begin{equation*}
    0\to H^1(G,\Fp)/\ann(a) \to H^2(G,\Fp)\to H^2(G,\Fp)/(a)\cdot
    H^1(G,\Fp)\to 0.
\end{equation*}
Assume that $d=\dim_{\Fp} H^2(G,\Fp) < \infty$, and set
\begin{equation*}
    x=\dim_{\Fp} H^1(G,\Fp)/\ann(a), \quad y=\dim_{\Fp}
    H^2(G,\Fp)/(a)\cdot H^1(G,\Fp).
\end{equation*} Then $d=x+y$ and we
have the following corollary on the size of $H^2(N,\Fp)$:
\begin{corollary*}
    Assume that $G$ and $N$ are as above.  Then
    \begin{equation*}
        \dim_{\Fp} H^2(N,\Fp) = x+ py.
    \end{equation*}
\end{corollary*}

Before the proof we need several intermediate results.  We assume
throughout this section that $F$ is a field containing a primitive
$p$-th root of unity $\xi_p$, $G=\Gal(F(p)/F)$, $N$ is a subgroup of
$G$ of index $p$ with fixed field $K=F(\root{p}\of{a})$, and
$\sigma$ denotes a fixed generator of $G/N$ with $\root{p}
\of{a}^{\sigma-1}=\xi_p$. For a field $F$, let $G_F$ denote its
absolute Galois group.  Observe that because
$1+\sigma+\cdots+\sigma^{p-1}\equiv (\sigma-1)^{p-1}$ modulo $p$,
the endomorphism $(\sigma-1)^{p-1}$ on $H^i(N,\Fp)$ is identical to
the composition $\res\circ\cor$.

\begin{proposition}\label{pr:kers}\
    \begin{enumerate}
        \item The inflation maps $\inf:H^i(G,\Fp)\to H^i(G_F,\Fp)$
        and $\inf:H^i(N,\Fp)\to H^i(G_K,\Fp)$, $i=1, 2$, are
        isomorphisms.  Moreover, the latter isomorphisms are
        $\Fp[G/N]$-equivariant.
        \item The kernel of the corestriction map $\cor:
        H^2(N,\Fp) \to H^2(G,\Fp)$ is $(\sigma-1)H^2(N,\Fp) + \res
        H^2(G,\Fp)$.
        \item The kernel of the restriction map $\res:
        H^2(G,\Fp)\to H^2(N,\Fp)$ is $(a) \cdot H^1(G,\Fp)$.
    \end{enumerate}
\end{proposition}

\begin{proof}
    (1). We prove first the statements for $G$ and $G_F$.  Observe
    that since $F$ contains a primitive $p$-th root of unity, $F(p)$
    is closed under taking $p$-th roots and hence
    $H^1(G_{F(p)},\Fp)= \{0\}$. Therefore by
    \cite[Theorem~11.5]{MeSu} we see that $H^2(G_{F(p)}, \Fp) =
    \{0\}$ as well. Then, considering the Lyndon-Hochschild-Serre
    spectral sequence associated to $1 \to G_{F(p)} \to G_F \to G
    \to 1$, we obtain that $\inf: H^i(G,\Fp) \to H^i(G_F,\Fp)$ is an
    isomorphism for each $i=1,2$.  The proof that
    $\inf:H^i(N,\Fp)\to H^i(G_K,\Fp)$, $i=1, 2$ are isomorphisms
    follows as above.  The fact that these isomorphisms are
    $\Fp[G/N]$-equivariant follows immediately from the explicit
    action of $\Fp[G/N]$ on cochains.

    (2). By \cite[Proposition~15.1]{MeSu}, the kernel of the
    corestriction map $\cor:H^2(G_K,\Fp) \to H^2(G_F,\Fp)$ is
    $(\sigma-1)H^2(G_K,\Fp) + \res H^2(G_F,\Fp)$. Hence the second
    row is exact in the following commutative diagram. (Observe that
    $\sigma$ commutes with $\inf$ by (1), and the right-hand square
    commutes by \cite[Proposition~1.5.5ii]{NSW}.)
    \begin{equation*}
        \xymatrix{{H^2(N,\Fp) \oplus H^2(G,\Fp)}
        \ar[r]^{\qquad\quad\stackrel{(\sigma-1)}{\oplus\res}}
        \ar[d]^{\inf\oplus\inf} & H^2(N,\Fp) \ar[r]^{\cor}
        \ar[d]^{\inf} & H^2(G, \Fp) \ar[d]^{\inf} \\
        H^2(G_K,\Fp)\oplus H^2(G_F,\Fp)
        \ar[r]^{\qquad\quad\stackrel{(\sigma-1)}{\oplus\res}} &
        H^2(G_K,\Fp) \ar[r]^{\cor} & H^2(G_F,\Fp)}
    \end{equation*}
    The first row is therefore exact and we have our statement.

    (3). By \cite[Proposition~5]{Me} and \cite[Theorem~11.5]{MeSu},
    the kernel of the restriction map $\res:H^2(G_F,\Fp) \to
    H^2(G_K,\Fp)$ is $(a)\cdot H^1(G_F,\Fp)$.  A commutative
    diagram analogous to that of part (2) then gives our
    statement.
\end{proof}

\begin{corollary*}
    Suppose that the corestriction map $\cor:H^2(N,\Fp)\to
    H^2(G,\Fp)$ is surjective.  Then $\ker \cor =
    (\sigma-1) H^2(N,\Fp)$.
\end{corollary*}

\begin{proof}
    By part (2) above, it is sufficient to show that $\res
    H^2(G,\Fp)$ is a subset of $(\sigma-1)H^2(N,\Fp)$.  Let
    $\alpha\in H^2(G,\Fp)$. By hypothesis, there exists $\beta\in
    H^2(N,\Fp)$ such that $\cor \beta = \alpha$. Recalling that
    $\res \cor = (\sigma-1)^{p-1}$, we see that $\res \alpha =
    (\sigma-1)^{p-1} \beta\in (\sigma-1)H^2(N,\Fp)$.
\end{proof}

\begin{lemma}\label{le:x1}
    Suppose that the corestriction map $\cor:H^2(N,\Fp)\to
    H^2(G,\Fp)$ is surjective.  Then there exists a trivial
    $\Fp[G/N]$-submodule $X$ of $H^2(N,\Fp)$ such that
    \begin{equation*}
        \cor: X \to (a) \cdot H^1(G,\Fp)
    \end{equation*}
    is an isomorphism.  In fact,  $\cor(H^2(N,\Fp)^{G/N})=(a)\cdot
    H^1(G,\Fp).$
\end{lemma}

\begin{proof}
    Let $\Ic$ be an $\Fp$-basis for $(a)\cdot H^1(G,\Fp) \subset
    H^2(G,\Fp)$.  For each $(a)\cdot (f)\in \Ic$ we will define an
    element $x_f\in H^2(N,\Fp)$ such that $\cor x_f = (a) \cdot (f)$
    and $(\sigma-1) x_f = 0$.  Then the $\Fp$-span $X$ of $x_f$ will
    be our required module $X$.  If $p=2$, then we proceed as
    follows. By hypothesis there exists $x_f\in H^2(N,\F_2)$ such
    that $\cor x_f = (a) \cdot (f)$. Then
    \begin{equation*}
        (\sigma-1)x_f = (\sigma+1)x_f = \res \cor x_f =
        \res((a)\cdot (f)) = 0,
    \end{equation*}
    and hence $x_f\in H^2(N,\F_2)^{G/N}$.

    Now suppose that $p>2$.  If $\res((\xi_p)\cdot (f))=0$ then set
    $x_f = (\root{p}\of{a}) \cdot (f)$.  Observe that in this case
    $x_f\in H^2(N,\Fp)^{G/N}$ and by the projection formula
    \cite[Proposition~1.5.3iv]{NSW}, we have $\cor x_f = (a)\cdot
    (f)$. Otherwise, by hypothesis there exists $\alpha\in
    H^2(N,\Fp)$ such that $\cor\alpha = (\xi_p) \cdot (f)$.  Let
    $\beta = (\sigma-1)^{p-2}\alpha$.  From $(\sigma-1)^{p-1} =
    \res \cor$ we obtain $(\sigma-1)\beta = \res((\xi_p)\cdot (f))$.
    Now set $x_f := (\root{p}\of{a}) \cdot (f) - \beta$.  Then
    \begin{equation*}
        (\sigma-1)x_f = \res((\xi_p) \cdot
        (f))-\res((\xi_p) \cdot (f)) = 0,
    \end{equation*}
    so $x_f\in H^2(N,\Fp)^{G/N}$. Observe that since the
    corestriction commutes with $\sigma$
    \cite[Proposition~1.5.4]{NSW}, $\cor$ vanishes on the image of
    $\sigma-1$.  Hence $\cor\beta=0$. By the projection formula
    again, $\cor x_f = (a) \cdot (f)$.

    Letting $X$ be the $\Fp$-span of the elements $x_f$, we have
    the first statement of the lemma.

    For the second statement, let $\gamma\in H^2(N,\Fp)^{G/N}$.
    Then $\res\cor \gamma=(\sigma-1)^{p-1}\gamma=0.$
    By Proposition~\ref{pr:kers}, part (3),
    \begin{equation*}
        \cor\gamma\in \ker \res=(a)\cdot H^1(G,\Fp).
    \end{equation*}
    Therefore $\cor(H^2(N,\Fp)^{G/N})\subset (a)\cdot H^1(G,\Fp).$
    The reverse inclusion follows from the first statement.
\end{proof}

\begin{lemma}\label{le:fixednorm}
    Suppose that the corestriction map $\cor:H^2(N,\Fp)\to
    H^2(G,\Fp)$ is surjective.
    Then
    \begin{equation*}
        H^2(N,\Fp)^{G/N} \cap (\sigma-1) H^2(N,\Fp) =
        (\sigma-1)^{p-1} H^2(N,\Fp).
    \end{equation*}
\end{lemma}

\begin{proof}
    Since
    \begin{equation*}
        (\sigma-1)^{p-1} H^2(N,\Fp) \subset H^2 (N,\Fp)^{G/N}
        \cap (\sigma-1) H^2(N,\Fp),
    \end{equation*}
    it is sufficient to prove the reverse inclusion. If $p=2$ the
    reverse inclusion is true since $(\sigma-1) H^2 (N,\Fp) \subset
    H^2 (N,\Fp)^{G/N}$. Therefore assume that $p>2$.

    Let
    \begin{equation*}
        \gamma \in H^2(N,\Fp)^{G/N} \cap (\sigma-1) H^2 (N,\Fp).
    \end{equation*}
    Since $0 \in (\sigma-1)^{p-1} H^2 (N,\Fp)$ we also assume that
    $\gamma \ne 0$. Then $\gamma = (\sigma-1) \beta$ for some $\beta
    \in H^2 (N,\Fp)$. We shall show by induction on $j$, $2 \le j \le
    p$, that there exists $\beta_j \in H^2 (N,\Fp)$ such that
    \begin{equation*}
        (\sigma-1)^{j-1} \beta_j = \gamma.
    \end{equation*}
    Then for $\beta_p$ we shall have
    \begin{equation*}
        (\sigma-1)^{p-1} \beta_p = \gamma \in (\sigma-1)^{p-1}
        H^2 (N,\Fp),
    \end{equation*}
    which will prove our desired inclusion
    \begin{equation*}
        H^2 (N,\Fp)^{G/N} \cap (\sigma-1) H^2 (N,\Fp)
        \subset (\sigma-1)^{p-1} H^2 (N,\Fp).
    \end{equation*}

    If $j = 2$ we set $\beta_2 = \beta$. Assume now that $2 \le j-1
    < p$ and that $(\sigma-1)^{j-2} \beta_{j-1} = \gamma$ for some
    $\beta_{j-1} \in H^2 (N,\Fp)$. Consider $\delta = \cor
    \beta_{j-1}$. Since
    \begin{equation*}
        (\sigma-1)^{j-1} \beta_{j-1} = (\sigma-1) \gamma
        = 0
    \end{equation*}
    and $(\sigma-1)^{p-1} = \res \cor$, we obtain $\res\cor
    \beta_{j-1} = \res \delta = 0$.  By Proposition~\ref{pr:kers},
    part (3), $\delta = (a) \cdot (f)$ for $(f)\in H^1(G,\Fp)$.  By
    Lemma~\ref{le:x1} there exists an element $x\in
    H^2(N,\Fp)^{G/N}$ such that $\cor x = (a) \cdot (f)$.  Let
    $\beta'_{j-1} = \beta_{j-1} - x$.

    From $(\sigma-1) x = 0$ and $j > 2$ we obtain
    \begin{equation*}
        (\sigma-1)^{j-2} \beta'_{j-1} = (\sigma-1)^{j-2}
        \beta_{j-1} = \gamma.
    \end{equation*}
    Moreover $\cor \beta'_{j-1} = 0$. By the Corollary to
    Proposition~\ref{pr:kers}, there exists $\beta_j \in H^2
    (N,\Fp)$ such that $(\sigma-1) \beta_j = \beta'_{j-1}$ and hence
    \begin{equation*}
        (\sigma-1)^{j-1} \beta_j = (\sigma-1)^{j-2} \beta'_{j-1} =
        \gamma,
    \end{equation*}
    as desired.
\end{proof}

\begin{lemma}\label{le:freecyc}
    Let $H$ be a cyclic group of order $p$ generated by $\sigma$, and
    let $T$ be an $\Fp[H]$-module.  Suppose that $\alpha \in T$ and
    $(\sigma-1)^{p-1} \alpha \ne 0$.  Then the $\Fp[H]$-submodule
    $\langle \alpha \rangle$ of $T$ generated by $\alpha$ is a free
    $\Fp[H]$-module.
\end{lemma}

\begin{proof}
    Let $S = \Fp[H]$ and let $I$ be any nonzero ideal of $S$. Let $w
    \ne 0$ be in $I$. Write
    \begin{equation*}
        w = \sum_{i=k}^{p-1} c_i (\sigma-1)^i, \quad k \in \{ 0,1,
        \dots,p-1 \}, \ c_i \in \Fp, \ c_k \ne 0.
    \end{equation*}
    Then also $w(\sigma-1)^{p-1-k} = c_k (\sigma-1)^{p-1} \in I$,
    and hence $(\sigma-1)^{p-1} \in I$.

    Now consider $\ann_S (\alpha) = \{ s \in S \ \vert\ s \alpha = 0
    \}$. If $\ann_S (\alpha) \neq \{0\}$ then $(\sigma-1)^{p-1} \in
    \ann_S (\alpha)$, contradicting our hypothesis. Hence $\ann_S
    (\alpha) = \{0\}$ and we see that $\langle \alpha \rangle$ is a
    free $\Fp[H]$-submodule of $T$.
\end{proof}

\begin{proof}[Proof of Theorem~\ref{th:h2}]
    By Lemma~\ref{le:x1}, there exists a trivial $\Fp[G/N]$-sub\-module
    $X$ of $H^2(N,\Fp)$ such that $\cor:X\to (a)\cdot H^1(G,\Fp)$ is an
    isomorphism.  Hence $\dim_{\Fp} X = \dim_{\Fp} \ H^1(G,\Fp)/
    \ann(a)$.  (Recall that $\ann(a) = \{(b)\in H^1(G,\Fp) \ \vert\
    (a)\cdot(b)= 0\}$.)

    Furthermore, there exists a maximal free $\Fp[G/N]$-submodule
    $Y$ of $H^2(N,\Fp)$, as follows.  ($Y$ may be zero since we
    consider $\{0\}$ to be a free $\Fp[G/N]$-module.) First by
    \cite[\S III.1, Proposition~1.4]{La}, an $\Fp[G/N]$-module $M$
    is free precisely when $H^2(G/N,M)=\{0\}$.  Observe that the
    trace map $1+ \sigma + \cdots+\sigma^{p-1}=(\sigma-1)^{p-1}$ in
    $\Fp [G/N]$. Recall that for any $\Fp[G/N]$-module $M$ we have
    $H^2(G/N,M) = M^{G/N}/(\sigma-1)^{p-1}M$.  (See \cite[I.5]{La}.)
    Therefore $M$ is a free $\Fp[G/N]$-module if and only if
    $M^{G/N} = (\sigma-1)^{p-1}M$.  Let $\Sc$ denote the set of free
    $\Fp[G/N]$-submodules of $H^2(N,\Fp)$.  Suppose $\Tc$ is a
    totally ordered subset of $\Sc$, and let $W=\cup_{S\in \Tc} S$.
    Then $W$ is the inductive limit of $S\in \Tc$.  Thus we have:
    \begin{equation*}
        H^2(G/N,W) = H^2(G/N,\varinjlim_{S\in \Tc} S) =
        \varinjlim_{S\in \Tc} H^2(G/N,S) = \{0\}.
    \end{equation*}
    Hence $W$ is a free $\Fp[G/N]$-module.  By Zorn's Lemma, $\Sc$
    contains a maximal element $Y$. We then have $Y^{G/N} =
    (\sigma-1)^{p-1}Y$.  Since $\dim_{\Fp} \Fp[G/N]^{G/N} =
    \dim_{\Fp} \langle (\sigma-1)^{p-1} \rangle = 1$, we obtain
    \begin{equation*}
        \rank Y = \dim_{\Fp} Y^{G/N} = \dim_{\Fp} (\sigma-1)^{p-1}Y.
    \end{equation*}

    Because free $\Fp[G/N]$-modules are injective (see
    \cite[Theorem~11.2]{Ca}) we may write $H^2(N,\Fp) = Y \oplus R$
    for some $\Fp[G/N]$-submodule $R$ of $H^2(N,\Fp)$. We will show
    that $R\cong X$ as $\Fp[G/N]$-modules.

    We first show that $R$ is a trivial $\Fp[G/N]$-module.  If there
    exists $\alpha\in R$ with $(\sigma-1)^{p-1}\alpha\neq 0$, by
    Lemma~\ref{le:freecyc} we see that $Y\oplus \langle \alpha\rangle$
    is a larger free $\Fp[G/N]$-submodule, a contradiction.  We obtain
    $(\sigma-1)^{p-1}R=\{0\}$ and $(\sigma-1)^{p-1}Y =
    (\sigma-1)^{p-1}H^2(N,\Fp)$.

    Because $(\sigma-1)^{p-1}R=\{0\}$ there exists a minimal $0\le
    l\le p-1$ such that $(\sigma-1)^l R=\{0\}$. Suppose $l>1$.  Then
    \begin{equation*}
        \{0\}\ne (\sigma-1)^{l-1} R\subset H^2(N,\Fp)^{G/N}\cap
        (\sigma-1)H^2(N,\Fp).
    \end{equation*}
    By Lemma~\ref{le:fixednorm},
    \begin{equation*}
        \{0\}\ne (\sigma-1)^{l-1} R\subset
        (\sigma-1)^{p-1}H^2(N,\Fp)=(\sigma-1)^{p-1}Y.
    \end{equation*}
    But then $\{0\}\ne (\sigma-1)^{l-1}R\subset R\cap Y$, a
    contradiction. Therefore $l\le 1$ and $(\sigma-1)R=\{0\}$.  Hence
    $R$ is indeed a trivial $\Fp[G/N]$-module.

    In fact, we claim that $R\cap (\sigma-1)H^2(N,\Fp)=\{0\}$. We
    have
    \begin{eqnarray*}
        R\cap (\sigma-1)H^2(N,\Fp)&\subset& H^2(N,\Fp)^{G/N}\cap
        (\sigma-1)H^2(N,\Fp)\\ &=&(\sigma-1)^{p-1}H^2(N,\Fp)=
        (\sigma-1)^{p-1}Y.
    \end{eqnarray*}
    From $R\cap Y=\{0\}$ we obtain $R\cap
    (\sigma-1)H^2(N,\Fp)=\{0\}$.

    Now consider the image of $\cor$ on
    \begin{equation*}
        H^2(N,\Fp)^{G/N}=R\oplus Y^{G/N}=R\oplus
        (\sigma-1)^{p-1}H^2(N,\Fp).
    \end{equation*}
    Observe that since the corestriction commutes with $\sigma$
    \cite[Proposition~1.5.4]{NSW}, $\cor$ vanishes on the image of
    $\sigma-1$.  By Lemma~\ref{le:x1}, we find that $\cor R =
    (a)\cdot H^1(G,\Fp) = \cor X$. But by the Corollary to
    Proposition~\ref{pr:kers} and the fact that $R\cap
    (\sigma-1)H^2(N,\Fp)=\{0\}$, we deduce that $\cor$ acts
    injectively on $R$.  Since, by Lemma~\ref{le:x1}, $\cor$ also
    acts injectively on $X$, we have that $R\cong X$. Hence we
    obtain that $H^2(N,\Fp)\cong X\oplus Y$.

    Now we determine the rank of $Y$. We have
    $(\sigma-1)^{p-1}H^2(N,\Fp) = (\sigma-1)^{p-1}Y$, and hence
    \begin{equation*}
        \rank Y = \dim_{\Fp} (\sigma-1)^{p-1}H^2(N,\Fp).
    \end{equation*}
    Using the hypothesis $\cor H^2(N,\Fp) = H^2(G,\Fp)$ together
    with $\res\cor = (\sigma-1)^{p-1}$, we obtain that
    $(\sigma-1)^{p-1}H^2(N,\Fp) = \res H^2(G,\Fp)$.  By
    Proposition~\ref{pr:kers}, part (3), the kernel of $\res$ is
    $(a)\cdot H^1(G,\Fp)$.  We deduce then that $\rank Y =
    \dim_{\Fp} H^2(G,\Fp)/(a)\cdot H^1(G,\Fp)$.
\end{proof}

\begin{theorem}\label{th:h2bis}
    Let $F$ be a field containing a primitive $p$-th root of unity,
    and suppose that $G=\Gal(F(p)/F)$.  Let $N$ be a subgroup of $G$
    of index $p$, and suppose that the corestriction map
    $\cor:H^2(N,\Fp)\to H^2(G,\Fp)$ is surjective.

    Suppose that $X$ and $Y$ are $\Fp[G/N]$-submodules of
    $H^2(N,\Fp)$ such that $X$ trivial and $Y$ is free.  Then $\cor
    X\subset (a)\cdot H^1(G,\Fp)$ and the following are equivalent:
    \begin{enumerate}
        \item $\cor:X\to (a)\cdot H^1(G,\Fp)$ is an isomorphism, and
        $Y$ is a maximal free submodule
        \item $H^2(N,\Fp)=X\oplus Y$.
    \end{enumerate}
\end{theorem}

\begin{proof}
    Since $X$ is a trivial $\Fp[G/N]$-module, $\res \cor X =
    (\sigma-1)^{p-1}X = \{0\}$.  By Proposition~\ref{pr:kers}, part
    (3), $\cor X\subset (a)\cdot H^1(G,\Fp)$.

    (1)$\implies$(2).  Suppose $w\in X\cap Y$. Since $X$ is a
    trivial $\Fp[G/N]$-module, $w\in Y^{G/N}$.  Then because $Y$ is
    a free $\Fp[G/N]$-module, $Y^{G/N}=(\sigma-1)^{p-1}Y$. In
    particular, $w\in (\sigma-1)Y$.  Since $\cor$ vanishes on the
    image of $\sigma-1$, $\cor w=0$, and because $\cor$ is injective on
    $X$, $w=0$.  Hence the submodule of $H^2(G,\Fp)$ generated by
    $X$ and $Y$ is $X \oplus Y$.

    Let $R$ be a trivial $\Fp[G/N]$-submodule of $H^2(N,\Fp)$ such
    that $\cor R=\cor X$ and $H^2(N,\Fp)=R\oplus Y$, as in the proof
    of Theorem~\ref{th:h2}.  Since $(\sigma-1)R=\{0\}$ we deduce that
    $(\sigma-1)^{p-1}Y=(\sigma-1)^{p-1}H^2(N,\Fp)$.

    To prove that $X \oplus Y = H^2(N,\Fp)$ it suffices to prove
    that $R \subset X \oplus Y$.  Let $r \in R$. Then there exists
    $x \in X$ such that $\cor r = \cor x$.  Thus $u = r - x \in H^2
    (N,\Fp)^{G/N}$ and $\cor u = 0$. By the Corollary to
    Proposition~\ref{pr:kers} we obtain that $u \in (\sigma-1)
    H^2(N,\Fp)$. Thus
    \begin{equation*}
        u \in H^2(N,\Fp)^{G/N} \cap (\sigma-1) H^2(N,\Fp),
    \end{equation*}
    and so by Lemma~\ref{le:fixednorm},
    \begin{equation*}
        u \in (\sigma-1)^{p-1} H^2(N,\Fp) = (\sigma-1)^{p-1}Y.
    \end{equation*}
    Hence $r \in X \oplus Y$ as required and we have $X\oplus Y
    = H^2(N,\Fp)$.

    (2)$\implies$(1).  By Lemma~\ref{le:x1}, $\cor
    (H^2(N,\Fp)^{G/N}) = (a)\cdot H^1(G,\Fp)$.  Since $Y$ is free,
    $Y^{G/N} = (\sigma-1)^{p-1}Y$, and since $\cor$ vanishes on the
    image of $\sigma-1$, $\cor Y^{G/N}=\{0\}$.  From
    $H^2(N,\Fp)^{G/N}=X\oplus Y^{G/N}$ we deduce that $\cor:X\to
    (a)\cdot H^1(G,\Fp)$ is surjective.  Now if $x\in X$ with $\cor
    x=0$ then by the Corollary to Proposition~\ref{pr:kers}, $x\in
    (\sigma-1)H^2(N,\Fp)$.  Because $X$ is trivial and $X\oplus
    Y=H^2(N,\Fp)$, we see that
    \begin{align*}
        x &\in (\sigma-1)H^2(N,\Fp)\cap H^2(N,\Fp)^{G/N}\\
        &=(\sigma-1)^{p-1}H^2(N,\Fp) = (\sigma-1)^{p-1}Y
    \end{align*}
    by Lemma~\ref{le:fixednorm}.  Then $x\in X\cap Y$, and so
    $x=0$. Hence $\cor$ is injective on $X$ and therefore $\cor:X\to
    (a)\cdot H^1(G,\Fp)$ is an isomorphism.

    Finally we show that $Y$ is a maximal free $\Fp[G/N]$-submodule.
    Suppose $Y \subset T$ where $T$ is a free $\Fp[G/N]$-submodule
    of $H^2(N,\Fp)$. Then because $Y$ is injective we can write $T =
    Y \oplus S$ for some $\Fp[G/N]$-module $S$. Then $S$ is a
    projective $\Fp[G/N]$-module, and since each projective
    $\Fp[G/N]$-module is free (see \cite[Proof of Theorem~11.2,
    pp.~70--71]{Ca}) we see that $S$ is in fact a free
    $\Fp[G/N]$-submodule of $T$. Then we have
    \begin{align*}
        \res \cor T &= \res \cor Y \oplus \res \cor S \\
        &= (\sigma-1)^{p-1} Y \oplus (\sigma-1)^{p-1} S.
    \end{align*}
    But since $H^2(N,\Fp) = X \oplus Y$ and $X$ is a trivial
    $\Fp[G/N]$-submodule of $H^2(N,\Fp)$ we see that
    \begin{align*}
        \res \cor H^2(N,\Fp) &= (\sigma-1)^{p-1} H^2 (N,\Fp) \\
        &= (\sigma-1)^{p-1} Y.
    \end{align*}
    Hence $(\sigma-1)^{p-1} S = \{0\}$.  Since $S$ is free, $S =
    \{0\}$.  Thus $Y$ is indeed a maximal free $\Fp[G/N]$-submodule
    of $H^2(N,\Fp)$.
\end{proof}

\section{Proof of Theorem~\ref{th:new}}\label{se:proofs}

    Let $N$ be a subgroup of $G$ of index $p$. Since $G$ has
    cohomological dimension $2$, the corestriction map $\cor:
    H^2(N,\Fp) \to H^2(G,\Fp)$ is surjective
    \cite[Proposition~3.3.8]{NSW}. By Theorem~\ref{th:h2} we have a
    decomposition $H^2(N,\Fp)=X\oplus Y$, where $X$ is a trivial
    $\Fp[G/N]$-module and $Y$ is a free $\Fp[G/N]$-module. Hence
    $H^2(N,\Fp)$ is trivial if and only if $H^2(N,\Fp)$ contains no
    nonzero free submodule.  We have established the first
    equivalence of the theorem.

    For the next assertion, observe that if $G$ is a Demu\v{s}kin
    group of cohomological dimension $2$ and $N$ is a subgroup of
    $G$ of index $p$, by \cite[Theorem~1]{DuLa}, the
    $\Fp[G/N]$-module $H^2(N,\Fp)$ is a trivial $\Fp[G/N]$-module.

    Conversely, by the definition of a Demu\v{s}kin group, it
    suffices to show that $\dim_{\Fp} H^2(G,\Fp)=1$ and $\gamma_F$
    is nondegenerate.  Consider the decomposition $H^2(N,\Fp)$
    obtained above, for $N$ an arbitrary subgroup of index $p$. Let
    $a\in F^\times$ be chosen so that the fixed field of $N$ is
    $F(\root{p}\of{a})$. Since we are assuming that $H^2(N,\Fp)$
    contains no nonzero free summand, from Theorem~\ref{th:h2} we
    obtain $\dim_{\Fp}H^2(G,\Fp)/(a)\cdot H^1(G,\Fp) = 0$, or
    $(a)\cdot H^1(G,\Fp) = H^2(G,\Fp)$.  Hence $\gamma_F$ is
    strongly regular. Moreover, $H^2(N,\Fp)$ has $\Fp$-dimension
    \begin{align*}
        \dim_{\Fp}H^1(G,\Fp)/\ann(a)&=\dim_{\Fp} \left((a) \cdot
        H^1(G,\Fp)\right)\\ &=\dim_ {\Fp}H^2 (G,\Fp).
    \end{align*}
    Suppose that the pairing $\gamma_F$ is degenerate.  Then for
    some nonzero $(a)\in H^1(G,\Fp)$ we have $(a)\cdot
    H^1(G,\Fp)=\{0\}$. Then $H^2(G,\Fp)=\{0\}$, contradicting the
    cohomological dimension of $G$.  Hence $\gamma_F$ is
    nondegenerate. Now if $p=2$ we have $H^2(G,\F_2)\cong \F_2$ by
    Proposition~\ref{pr:p2}. If $p>2$ and we assume the Elementary
    Type Conjecture, then by Theorem~\ref{th:notstrreg}, $\gamma_F$
    is of $p$-local type and hence $\dim_{\Fp} H^2(G,\Fp)=1$.

    Thus $G$ is a Demu\v{s}kin group as required. \qed

\begin{remark}
    In the proof of Theorem~\ref{th:new} we cited \cite[Theorem
    1]{DuLa} to establish that if $G$ is Demu\v{s}kin with
    $\cd(G)=2$, then $H^2(N,\Fp)$ is a trivial $\Fp[G/N]$-module.
    This result also follows from the fact that open subgroups of
    Demu\v{s}kin groups $G\neq \Z/2\Z$ are also Demu\v{s}kin
    \cite[Corollary I.4.5]{S2}. We observe that we can also deduce
    this result in our setting when $G = \Gal(F(p)/F)$ from
    Theorem~\ref{th:h2}, as follows. By Theorem~\ref{th:h2},
    $H^2(N,\Fp)$ is the direct sum of a trivial $\Fp[G/N]$-module
    $X$ and a free $\Fp[G/N]$-module $Y$. Since $G$ is Demu\v{s}kin,
    $\gamma_F$ is strongly regular.  From Theorem~\ref{th:h2}(2), we
    have $Y = \{0\}$. Hence $H^2(N,\Fp)$ is a trivial
    $\Fp[G/N]$-module as required. More precisely, from
    Theorem~\ref{th:h2}(1) and the fact that $\gamma_F$ is strongly
    regular we obtain $H^2(N,\Fp)\cong X\cong \Fp$.
\end{remark}

\begin{remark}\label{re:etc}
    By Theorem~\ref{th:notstrreg}, the Elementary Type Conjecture
    for Odd $p$ holds for a field $F$ with a strongly regular
    non-totally degenerate $p$-quaternionic pairing $\gamma_F$ if
    and only if $G=\Gal(F(p)/F)$ is Demu\v{s}kin. Thus
    Theorem~\ref{th:new} may be viewed as a translation of the
    Elementary Type Conjecture to the language of Galois
    $\Fp[G/N]$-modules $H^2(N,\Fp)$ in the case of strongly regular
    non-totally degenerate $p$-quaternionic pairings. There is some
    additional interest in this formulation because $p$-quaternionic
    pairings which are strongly regular but not weakly $p$-local
    have been abstractly constructed (see \cite[Theorem~9]{Ku}), and
    it is not known whether these pairings are realizable as
    $\gamma_F$ for suitable fields $F$.
\end{remark}

\section{Structure of $H^1(N,\Fp)$}

In this section we keep our assumption that a primitive $p$-th root
of unity lies in $F$. For any finitely generated pro-$p$-group $T$
we set $d(T) = \dim_{\Fp} H^1 (T,\Fp)$.

If $G$ is a Demu\v{s}kin pro-$p$-group then it is well-known
that
  \begin{equation*}
    d(N) = p (d(G) - 2) + 2
  \end{equation*}
for any subgroup $N$ of index $p$ of $G$. Moreover this formula
characterizes Demu\v{s}kin groups among finitely generated
pro-$p$-groups $G$ with $\dim_{\Fp} H^2 (G,\Fp) = 1$. (See
\cite{DuLa} or \cite[Theorem~3.9.15]{NSW}.) In this section we show
that this formula has an attractive explanation when $G =
\Gal(F(p)/F)$. In the following theorem $K$ is the fixed field in
$F(p)$ of the index $p$ subgroup $N$ of $G$.

\begin{theorem}\label{th:h1}
    Let $F$ be a field containing a primitive $p$-th root of unity
    $\xi_p$, and suppose that $G=\Gal(F(p)/F)$ is a Demu\v{s}kin
    group of rank $d(G) = \dim_{\Fp} H^1(G,\Fp) = n$.

    If $p>2$, then for each subgroup $N$ of $G$ of index $p$ we
    have a decomposition into $\Fp[G/N]$-modules
    \begin{equation*}
        H^1(N,\Fp) = X \oplus Y
    \end{equation*}
    where $X$ is an $\Fp[G/N]$-module of dimension 2 and $Y$ is
    a free $\Fp[G/N]$-module of rank $n-2$.  The module $X$ is
    trivial if $\xi_p\in N_{K/F}(K^\times)$ and is cyclic of
    dimension $2$ otherwise.

    If $p=2$ then for each subgroup $N$ of $G$ of index $p$ we
    have one of two decompositions into $\F_2[G/N]$-modules
    \begin{equation*}
        H^1(N,\F_2) = X \oplus Y \text{\ \ or\ \ } H^1(N,\F_2) = Y.
    \end{equation*}
    The first case occurs when $-1\in N_{K/F}(K^\times)$, and then
    $X$ is trivial of dimension 2 and $Y$ is free of rank $n-2$.
    The second occurs when $-1\not\in N_{K/F}(K^\times)$, and then
    $Y$ is free of rank $n-1$.
\end{theorem}

\begin{proof}
    Observe that for $N$ an index $p$ subgroup of the Demu\v{s}kin
    group $G$ and $K$ its fixed field in $F(p)$, we have $\dim_{\Fp}
    F^\times/ N_{K/F}(K^\times) = 1$.  Using equivariant Kummer
    theory, as explained in \cite{W2}, to identify the first
    cohomology groups with their corresponding $p$th-power classes
    as $\Fp[G/N]$-modules, the result then follows from the
    determination of the $\Fp[G/N]$-module structure of
    $K^\times/K^{\times p}$ in \cite[Theorem~3]{MiS}.
\end{proof}

\section*{Acknowledgements}

We thank Professors I.~Efrat, M.~Kula, J.-P.~Serre, and the referee
for their careful reading of preliminary versions of this paper and
for their valuable comments and suggestions. We are also grateful to
M.~Fried for his encouragement and infectious enthusiasm regarding
this research.

\end{document}